\documentclass[12pt,reqno]{amsart}

\DeclareMathOperator{\ann}{ann}%
\DeclareMathOperator{\cd}{cd}%
\DeclareMathOperator{\cf}{cf}%
\DeclareMathOperator{\ct}{ct}%
\DeclareMathOperator{\chr}{char}%
\DeclareMathOperator{\cor}{cor}%
\DeclareMathOperator{\res}{res}%
\DeclareMathOperator{\supp}{supp}%

\usepackage{amssymb}
\usepackage{hyperref}

\newcommand{\C}{\mathbb{C}}
\newcommand{\F}{\mathbb{F}}
\newcommand{\Fp}{\F_p}
\newcommand{\Gal}{\text{\rm Gal}}
\newcommand{\N}{\mathbb{N}}
\newcommand{\Q}{\mathbb{Q}}
\newcommand{\Z}{\mathbb{Z}}

\begin{document}

\title{When is Galois Cohomology Free or Trivial?}

\author[Nicole Lemire]{Nicole Lemire$^\star$}
\address{Department of Mathematics, Middlesex College, \
University of Western Ontario, London, Ontario \ N6A 5B7 \ CANADA}
\thanks{$^\star$Research supported in part by NSERC grant R3276A01.}
\email{nlemire@uwo.ca}

\author[J\'{a}n Min\'{a}\v{c}]{J\'an Min\'a\v{c}$^{*\dagger}$}
\address{Department of Mathematics, Middlesex College, \
University of Western Ontario, London, Ontario \ N6A 5B7 \ CANADA}
\thanks{$^*$Research supported in part by NSERC grant R0370A01 and
by a Distinguished Professorship during 2004--2005 at the University
of Western Ontario.}
\thanks{$^\dag$Supported by the Mathematical Sciences Research
Institute, Berkeley.} \email{minac@uwo.ca}

\author[John Swallow]{John Swallow$^\ddag$}
\address{Department of Mathematics, Davidson College, Box 7046,
Davidson, North Carolina \ 28035-7046 \ USA}
\thanks{$^\ddag$Research supported in part by NSA grant
MDA904-02-1-0061.} \email{joswallow@davidson.edu}

\begin{abstract}
Let $p$ be a prime and $F$ a field containing a primitive $p$th root
of unity.  Let $E/F$ be a cyclic extension of degree $p$ and $G_E
\triangleleft G_F$ the associated absolute Galois groups. We
determine precise conditions for the cohomology group
$H^n(E)=H^n(G_E,\Fp)$ to be free or trivial as an
$\Fp[\Gal(E/F)]$-module.  We examine when these properties for
$H^n(E)$ are inherited by $H^k(E)$, $k>n$, and, by analogy with
cohomological dimension, we introduce notions of cohomological
freeness and cohomological triviality. We give examples of $H^n(E)$
free or trivial for each $n\in \N$ with prescribed cohomological
dimension.
\end{abstract}

\date{February 25, 2005}

\maketitle

\newtheorem{theorem}{Theorem}
\newtheorem{proposition}{Proposition}
\newtheorem{lemma}{Lemma}
\newtheorem{corollary}{Corollary}

\theoremstyle{definition}
\newtheorem*{remark*}{Remark}
\newtheorem*{example*}{Example}

\parskip=9pt plus 2pt minus 2pt

Let $p$ be a prime and $F$ a field containing a primitive $p$th root
of unity $\xi_p$. Let $E/F$ be a cyclic extension of degree $p$ and
$G_E$ the absolute Galois group of $E$.  In our previous paper
\cite{LMS} we determined the structure of $H^n(G_E, \Fp)$, $n\in
\N$, as an $\Fp[G]$-module.  In this paper we study more closely the
question of when $H^n(G_E, \Fp)$ is free or trivial as an
$\Fp[G]$-module.

Let $a\in F$ satisfy $E=F(\root{p}\of{a})$. We write $H^n(F)$ for
$H^n(G_F, \Fp)$ and $\ann_n x$ for the annihilator of $x$ under
the cup-product operation on $H^n(F)$. (Thus $\ann_n x \subset
H^n(F)$.) Let $(f)\in H^1(F)$ denote the class of $f$ under the
Kummer isomorphism of $H^1(F)$ with the $p$th-power classes of
$F^\times:= F \setminus \{0\}$, and let $(f,g) \in H^2(F)$ denote
the cup-product of $(f)$ and $(g)\in H^1(F)$.

We first give precise conditions for free $\Fp[G]$-module
cohomology.

\begin{theorem}\label{th:free}  Let $n\in \N$.

    Suppose $p>2$.  Then the following are equivalent:
    \begin{enumerate}
        \item $H^n(E)$ is a free $\Fp[G]$-module
        \item $H^{n-1}(F) = \ann_{n-1} (a)$
        \item $\res\colon H^n(F)\to H^n(E)$ is injective
        \item $\cor\colon H^{n-1}(E)\to H^{n-1}(F)$ is
        surjective.
    \end{enumerate}

    Suppose $p=2$.  Then the following are equivalent:
    \begin{enumerate}
        \item $H^n(E)$ is a free $\F_2[G]$-module
        \item $\ann_{n-1} (a) = \ann_{n-1} (a,-1)$ and \\
        $H^n(F) = \cor H^n(E) + (a) \cup H^{n-1}(F)$
        \item $\ann_{n-1} (a) = \ann_{n-1} (a,-1)$ and \\
        $H^n(F)= \ann_n (a) + (a) \cup H^{n-1}(F)$
        \item $H^n(F) = \ann_n (a) \oplus (a)\cup H^{n-1}(F)$.
    \end{enumerate}
\end{theorem}

In the following theorem we examine to what extent free cohomology
is hereditary.

\begin{theorem}\label{th:freeher}
    Suppose that either
    \begin{itemize}
        \item $p>2$ or
        \item $p=2$ and $a\in (F^{\times 2}+F^{\times 2})
        \setminus F^2$.
    \end{itemize}

    Then free cohomology is hereditary:  if $n\in \N$, then
    for all $m\ge n$,
    \begin{equation*}
        H^n(E) \text{ is a free $\Fp[G]$-module} \implies H^m(E)
        \text{ is a free $\Fp[G]$-module}.
    \end{equation*}

    Moreover, if $H^m(E)$, $m\in \N$, is a free $\Fp[G]$-module,
    then the sequence
    \begin{equation*}
        0\to H^m(F) \xrightarrow{\res} H^m(E)
        \xrightarrow{\cor} H^m(F) \to 0.
    \end{equation*}
    is exact in the first and third terms.
\end{theorem}

We consider Theorems~\ref{th:free} and \ref{th:freeher} in
section~\ref{se:free}.  We moreover show that when $p>2$, $H^1(E)$
is never free.  When $p=2$ we show that free cohomology is not
generally hereditary and establish a condition for hereditary
freeness that is more general than the one given above.

We next give precise conditions for trivial $\Fp[G]$-module
cohomology.

\begin{theorem}\label{th:trivial}  Let $n\in \N$.

    Suppose $p>2$.  Then the following are equivalent:
    \begin{enumerate}
        \item $H^n(E)$ is a trivial $\Fp[G]$-module
        \item $(\xi_p)\cup H^{n-1}(F) \subset (a) \cup H^{n-1}(F)$
        and \\ $\ann_n (a) = (a) \cup H^{n-1}(F)$
        \item $(\xi_p)\cup H^{n-1}(F) \subset (a) \cup H^{n-1}(F)$
        and \\ $H^n(E) = \res H^n(F) + (\root{p}\of{a}) \cup \res
        H^{n-1}(F)$.
    \end{enumerate}

    Suppose $p=2$.  Then the following are equivalent:
    \begin{enumerate}
        \item $H^n(E)$ is a trivial $\F_2[G]$-module
        \item $\ann_n (a) \subset (a) \cup \ann_{n-1} (a,-1)$.
    \end{enumerate}

    In the $p=2$ case, suppose additionally that $a\in (F^{\times
    2}+F^{\times 2}) \setminus F^2$. Then the conditions above are
    also equivalent to
    \begin{enumerate}
    \setcounter{enumi}{2}
        \item $H^n(E) = \res H^n(F) + (\delta) \cup \res
        H^{n-1}(F)$ where $(\delta)\in H^1(E)^G$ satisfies
        $N_{E/F}(\delta) = (a)$.
    \end{enumerate}
\end{theorem}

For $p>2$ and $n=1$ the second condition in (3) was observed in
\cite[Lemma~3]{War}.

We deduce that trivial $\Fp[G]$-module cohomology is a hereditary
property.

\begin{theorem}\label{th:trivher}
    Trivial $\Fp[G]$-module cohomology is hereditary: if $n\in \N$,
    then for all $m\ge n$,
    \begin{equation*}
        H^n(E)^G=H^n(E) \implies H^m(E)^G=H^m(E).
    \end{equation*}

    Moreover, if $H^m(E)^G=H^m(E)$, $m\in \N$, then the following
    sequence is exact:
    \begin{multline*}
        0 \to \ann_{m-1} (a) \to H^{m-1}(F)
        \xrightarrow{(a)\cup -\ } H^m(F) \xrightarrow{\res} \\
        H^m(E) \xrightarrow{\cor} (a) \cup \ann_{m-1}
        \left((a)\cup (\xi_p)\right) \to 0,
    \end{multline*}
    where the map $\ann_{m-1} (a) \to H^{m-1}(F)$ is the natural
    inclusion.
\end{theorem}

We consider Theorems~\ref{th:trivial} and \ref{th:trivher} in
section~\ref{se:trivial}.

In section~\ref{se:exfree} we introduce $\cf(E/F)$, the largest
degree $n\in \N$ for which $H^n(E)$ is not free or $\infty$ if
$H^n(E)$ is never free, and we give examples, for each $m \ge
n\ge 1$, of extensions $E/F$ with $\cf(E/F)=n$ and $G_E$ a
pro-$p$-group of cohomological dimension $m$.

In section~\ref{se:extriv} we introduce $\ct(E/F)$, the largest
degree $n\in \N$ for which $H^n(E)$ is not a trivial $\Fp[G]$-module
or $\infty$ if $H^n(E)$ is never trivial, and we give examples,
for each $m\ge n\ge 1$, of extensions $E/F$ with $\ct(E/F)=n$
and $G_E$ a pro-$p$-group of cohomological dimension $m$.

Our proof relies on two recent results of Voevodsky in his proof
of the Bloch-Kato Conjecture. (Before Voevodsky's proof these
results were standard conjectures in Galois cohomology and they
were proved in important special cases.) In section~\ref{se:bkktheory}
we recall these results and present two corollaries deducing
collections of equivalent statements in Milnor $K$-theory. In
section~\ref{se:nl} we introduce various lemmas that give
sufficient conditions for our $\Fp[G]$-modules to be free or
trivial, demonstrate that some properties in Milnor $K$-theory are
hereditary, and establish some basic facts about certain
$p$-henselian fields we will use to construct our examples in
sections~\ref{se:exfree} and \ref{se:extriv}. For the convenience
of the reader we have made our paper quite independent of [LMS].

\section{Bloch-Kato and Milnor $K$-theory}\label{se:bkktheory}

The main ingredient for our determination of the $G$-module
structure of $H^n(E)$ is Milnor $K$-theory. (See \cite{Mi} and
\cite[Chap.~IX]{FV}.)  For $i\ge 0$, let $K_iF$ denote the $i$th
Milnor $K$-group of the field $F$, with standard generators
denoted by $\{f_1,\dots,f_i\}$, $f_1, \dots, f_i\in F\setminus
\{0\}$. For $\alpha\in K_iF$, we denote by $\bar\alpha$ the class
of $\alpha$ modulo $p$, and we use the usual abbreviation $k_nF$
for $K_nF/pK_nF$. The image of an element $\alpha\in K_iF$ in
$H^i(F)$ we also denote by $\alpha$. Because we will often use the
elements $\overline{\{a\}}$, $\overline{\{\xi_p\}}$,
$\overline{\{a,a\}}$, and $\overline{\{a,\xi_p\}}$, we omit the
bars for these elements. We also omit the bar in the element
$\overline{\{\root{p}\of{a}\}}$.

We write $N_{E/F}$ for the norm map $K_nE\to K_nF$, and we use the
same notation for the induced map modulo $p$.  We write
$\cor=\cor_{E/F}$ for the corresponding map of cohomology
$H^n(E)\to H^n(F)$.  We denote by $i_E$ the natural homomorphism
from $K_nF$ to $K_nE$, and we use the same notation for the
induced map modulo $p$. We denote by $\res=\res_{E/F}$ the
corresponding map of cohomology $H^n(F)\to H^n(E)$. We use
a well-known projection formula in Milnor $K$-theory several
times. (See \cite[page~81]{FW}.)

Our proof relies on the following two results in Voevodsky's proof
of the Bloch-Kato Conjecture. The first is the Bloch-Kato Conjecture
itself and its closely related Hilbert Theorem 90 for $K_m$:

\begin{theorem}[{\cite[Lemma~6.11 and \S 7]{Vo1} and \cite[\S 6 and
Theorem~7.1]{Vo2}}]\label{th:bk}\
    \begin{enumerate}
    \item Let $F$ be a field of characteristic not $p$ and
    $m\in\N$.  Then the norm residue homomorphism
    \begin{equation*}
        k_m F\to H^m(G_F, \mu_p^{\otimes m})
    \end{equation*}
    is an isomorphism.

    \item For any cyclic extension $E/F$ of degree $p$, the
    sequence
    \begin{equation*}
        K_mE \xrightarrow{1-\sigma} K_mE \xrightarrow{N_{E/F}}
        K_mF
    \end{equation*}
    is exact.
    \end{enumerate}
\end{theorem}

\noindent The second result establishes an exact sequence
connecting $k_mF$ and $k_mE$ for consecutive $m$.  (We translate
the statement of the original result to $K$-theory using the
previous theorem.)  In the following result $a$ is chosen to
satisfy $E=F(\root{p}\of{a})$.

\begin{theorem}[{\cite[Definition~5.1 and
Proposition~5.2]{Vo1}}]\label{th:es}
    Let $F$ be a field of characteristic not $p$ with no
    extensions of degree prime to $p$.  Then for any cyclic
    extension $E/F$ of degree $p$ and $m\ge 1$, the sequence
    \begin{equation*}
        k_{m-1}E \xrightarrow{N_{E/F}} k_{m-1}F \xrightarrow{\{a\}
        \cdot -\ } k_m F \xrightarrow{i_E} k_m E
    \end{equation*}
    is exact.
\end{theorem}

Now we observe that we may remove the hypothesis that the field $F$
has no extensions of degree prime to $p$.

\begin{theorem}[{Modification of
Theorem~\ref{th:es}: \cite[Theorem~5]{LMS}}]\label{th:esext}
    Let $F$ be a field containing a primitive $p$th root of unity.
    Then for any cyclic extension $E/F$ of degree $p$ and $m\geq 1$
    the sequence
    \begin{equation*}
        k_{m-1}E \xrightarrow{N_{E/F}} k_{m-1}F \xrightarrow{\{a\}
        \cdot -\ } k_m F \xrightarrow{i_E} k_m E
    \end{equation*}
    is exact.
\end{theorem}

We have the following corollaries of Theorem~\ref{th:esext}. For
an element $\bar\alpha$ of $k_iF$, let
\begin{equation*}
    \ann_{n-1} \bar\alpha = \ann_{k_{n-1}F}\bar\alpha = \ann
    \left(k_{n-1}F \xrightarrow{\bar\alpha \cdot -\ }
    k_{n-1+i}F\right)
\end{equation*}
denote the annihilator of the product with $\bar\alpha$.

\begin{corollary}\label{co:injsur}
    Assume the same hypotheses.  The following are equivalent for
    $n\in \N$:
    \begin{enumerate}
        \item $k_{n-1}F = \ann_{n-1} \{a\}$
        \item $k_{n-1}F = \ann_{n-1} \{a\} = \ann_{n-1}
        \{a,\xi_p\}$
        \item $i_E\colon k_nF\to k_nE$ is injective
        \item $N_{E/F}\colon k_{n-1}E\to k_{n-1}F$ is surjective.
    \end{enumerate}
\end{corollary}

\begin{proof}
    The equivalence of the items (1), (3), and (4) follows
    directly from the exact sequence.  Assuming (1) we see that
    \begin{equation*}
        k_{n-1}F = \ann_{n-1} \{a\} \subset \ann_{n-1} \{a,
        \xi_p\} \subset k_{n-1}F,
    \end{equation*}
    whence (2) follows, and (2) implies (1) trivially.
\end{proof}

In Lemma~\ref{le:normher} we show that all of the properties in
Corollary~\ref{co:injsur} are hereditary.

\begin{corollary}\label{co:cop2}
    Assume the same hypotheses. The following are equivalent for
    $n\in \N$:
    \begin{enumerate}
        \item $\ann_{n-1} \{a\} = \ann_{n-1} \{a,-1\}$ and \\
        $k_nF = N_{E/F}k_nE + \{a\}\cdot k_{n-1}F$
        \item $\ann_{n-1} \{a\} = \ann_{n-1} \{a,-1\}$ and \\
        $k_nF = \ann_{n} \{a\} + \{a\}\cdot k_{n-1}F$
        \item $k_nF = \ann_n \{a\} \oplus \{a\}\cdot k_{n-1}F$.
    \end{enumerate}
\end{corollary}

\begin{proof}
    (1)$\implies$(2).  This implication follows directly from
    $N_{E/F}k_nE = \ann_n \{a\}$.

    (2)$\implies$(3). Let $\bar\alpha\in (\{a\} \cdot k_{n-1}F)
    \cap \ann_n \{a\}$. Then $\bar\alpha = \{a\} \cdot \bar f$ for
    some $f\in K_{n-1}F$.  Since $\{a\}\cdot \bar\alpha = 0$,
    $\{a,a\}\cdot \bar f = 0$. Because $\{a,a\} = \{a,-1\}$, we have
    $\{a,-1\}\cdot \bar f = 0$, and by the first hypothesis,
    $\{a\}\cdot \bar f = 0$. Then $\bar\alpha = 0$ and the sum is
    direct.

    (3)$\implies$(1).  The second claim follows from the fact that
    $\ann_n \{a\} = N_{E/F}k_nE$. For the first, suppose
    $\{a,-1\}\cdot \bar f = 0$ for $f\in K_{n-1}F$. Because
    $\{a,-1\} = \{a,a\}$, we have
    \begin{equation*}
        \{a\}\cdot \bar f \in \left(\ann_n \{a\}\right) \cap
        \left( \{a\}\cdot k_{n-1}F \right) = \{0\}.
    \end{equation*}
    Hence $\bar f\in \ann_{n-1} \{a\}$ and $\ann_{n-1} \{a\} =
    \ann_{n-1} \{a, -1\}$ as required.
\end{proof}

\section{Notation and Lemmas}\label{se:nl}

For a field $F$, we let $F^\times$ denote its multiplicative group
$F\setminus\{0\}$. For the remainder of the paper $n\in \N$ denotes
an arbitrary natural number, $E/F$ a cyclic extension of fields of
degree $p$ with a primitive $p$th root of unity $\xi_p\in F$, and
$a\in F^\times$ an element such that $E=F(\root{p}\of{a})$.  Let
$G=\Gal(E/F)$, and choose $\sigma\in G$ to satisfy
$\root{p}\of{a}^{\sigma-1} = \xi_p$.  For $f,g\in F^\times$, we
write $(f)$ for the class of $f$ in $H^1(F)\cong F^\times/F^{\times
p}$ and $(f,g)$ for $(f)\cup (g) \in H^2(F)$.

\subsection{Module Structure}\

For $\gamma\in K_nE$, let $l(\gamma)$ denote the dimension of the
cyclic $\Fp[G]$-submodule $\langle \bar \gamma \rangle$ of $k_n E$
generated by $\bar\gamma$. Then we have, for $l(\gamma)\ge 1$,
\begin{equation*}
    (\sigma-1)^{l(\gamma)-1}\langle \bar\gamma \rangle = \langle
    \bar\gamma\rangle^G \neq 0 \text{\ \ \ and \ \ }
    (\sigma-1)^{l(\gamma)}\langle \bar\gamma \rangle=0.
\end{equation*}

We denote by $N$ the map $(\sigma-1)^{p-1}$ on $k_n E$. Because
$(\sigma-1)^{p-1}=1+\sigma+\dots+\sigma^{p-1}$ in $\Fp[G]$, we may
use $i_EN_{E/F}$ and $N$ interchangeably on $k_n E$.

Our first lemma establishes that in certain situations, all
elements in $(k_nE)^G$ are norm classes.

\begin{lemma}\label{le:extend}
    Let $n\in \N$.  Suppose that either
    \begin{itemize}
        \item $p>2$ and $N_{E/F}\colon k_{n-1}E\to k_{n-1}F$ is
        surjective, or
        \item $p=2$, $\ann_{n-1} \{a\} = \ann_{n-1} \{a,-1\}$, and
        \\ $k_nF = N_{E/F}k_nE + \{a\}\cdot k_{n-1}F$.
    \end{itemize}
    Then we have
    \begin{enumerate}
        \item For each $\gamma\in K_n E$, there exists $\alpha \in
        K_nE$ such that
        \begin{equation*}
            \langle N\bar\alpha \rangle = \langle \bar
            \gamma \rangle^G.
        \end{equation*}
        \item $(k_nE)^G = i_EN_{E/F}k_nE = (\sigma-1)^{p-1}k_nE =
        i_Ek_nF$.
    \end{enumerate}
\end{lemma}

\begin{proof}
    (1). Assume first that $p>2$.  By hypothesis, $N_{E/F} \colon
    k_{n-1}E \to k_{n-1}F$ is surjective, and then using the
    projection formula (\cite[p.~81]{FW}) we see that
    $N_{E/F}\colon k_nE \to k_nF$ is also surjective. Hence if
    $\bar\gamma\in i_Ek_nF$ then there exists $\bar\alpha\in k_nE$
    such that $N\bar\alpha = \bar\gamma$ and we are done.
    Otherwise, let $l=l(\gamma)$ and suppose $\bar\gamma\not\in
    i_Ek_nF$ and $1\le l\le i\le p$.

    If $l\ge 2$ we show by induction on $i$ that there exists
    $\alpha_i\in K_nE$ such that $\langle(\sigma-1)^{i-1}
    \bar\alpha_i\rangle = \langle \bar\gamma \rangle^G$.  Then
    setting $\alpha:=\alpha_p$, the proof will be complete in the
    case when $2\le l$.  The case $l=1$ we then handle at the end of
    the proof, using the case $2\le l$.

    Assume then that $l\ge 2$.  If $i=l$ then $\alpha_i= \gamma$
    suffices. Assume now that $1\le l \le i<p$ and that our
    statement is true for $i$. Set $c=N_{E/F}\alpha_i$. Since $i_E
    \bar c = N\bar\alpha_i = (\sigma-1)^{p-1} \bar\alpha_i$ and
    $i<p$, $i_E\bar c = 0$.

    By Corollary~\ref{co:injsur}, we have $\bar c = 0$, that is,
    $c = pf$ for some $f\in K_nF$.  Hence
    \begin{equation*}
        N_{E/F}\big(\alpha_i-i_E(f) \big) = 0.
    \end{equation*}
    By Theorem~\ref{th:bk}, there exists $\omega\in K_nE$ such
    that
    \begin{equation*}
        (\sigma-1)\omega = \alpha_i - i_E(f).
    \end{equation*}
    If $i>1$, $l(\alpha_i)>1$ and so $l(\alpha_i-i_E(f))>1$.
    Hence $(\sigma-1)^2\bar \omega = (\sigma-1)\bar \alpha_i \neq
    0$. Therefore $\langle (\sigma-1)^i\bar \omega \rangle = \langle
    \bar\gamma \rangle^G$ and we can set $\alpha_{i+1}=\omega$.

    Therefore we have proved that if $l(\gamma)\ge 2$ then there
    exists $\alpha\in K_n E$ such that $N\bar\alpha=(\sigma-1)^
    {l(\gamma)-1}\bar\gamma$.

    Now assume that $l(\gamma)=1$ but $\bar\gamma\notin i_E k_n F$.
    Then $\bar\gamma=\bar\alpha_1$ and $(\sigma-1)\bar\omega=
    \bar\alpha_1-i_E(\bar f)\neq 0$. Thus $l(\omega)=2$ and our
    argument above shows that there exists $\beta\in K_n E$ such
    that $N\bar\beta=(\sigma-1)\bar\omega=\bar\alpha_1-i_E(\bar
    f)$. As we observed at the beginning of our proof there exists
    an element $\delta\in K_n E$ such that $N\bar\delta=i_E(\bar f)$.
    Therefore we have:
    \begin{equation*}
        N(\bar\beta+\bar\delta)=\bar\alpha_1=\bar\gamma.
    \end{equation*}
    Thus we have established in all cases that for each $\gamma\in
    K_n E$ there exists $\alpha\in K_n E$ such that $\langle N\bar
    \alpha \rangle=\langle\bar\gamma \rangle^G$.

    Now consider the case $p=2$. In this case from our hypothesis
    $k_n F = N_{E/F}k_n E + \{a\}\cdot k_{n-1}F$ we again have $i_E
    N_{E/F} k_n E = i_E k_n F$. Therefore if $\bar\gamma \in i_E
    k_n F$ our statement follows. Assume that $\bar\gamma \in k_n
    E\setminus i_E k_n F$. Then $l(\gamma) \leq 2$, and if
    $l(\gamma)=2$ we may set $\alpha=\gamma$ and (1) follows
    again. Next we shall assume that $l(\gamma)=1$ and therefore
    $\bar\gamma\in (k_n E)^G$. Set $c = N_{E/F}\gamma$.  Then
    $i_E \bar c =0$.

    From Theorem~\ref{th:esext}, we conclude that $c=\{a\}\cdot
    g+2f$ for $g\in K_{n-1}F$ and $f\in K_nF$. Hence from the
    projection formula,
    \begin{align*}
        N_{E/F}\left(\gamma-\{\sqrt{a}\}\cdot i_E
        (g)-i_E(f)\right)  &=\left(\{a\}\cdot g+ 2f
        \right)-\{-a\}\cdot g-2f \\ &=\{-1\} \cdot g.
    \end{align*}
    Using Theorem~\ref{th:esext} again, we obtain that $\{a,-1\}
    \cdot \bar g = 0$. Our hypothesis $\ann_{n-1}\{a\} =
    \ann_{n-1}\{a,-1\}$ gives us that $\{a\} \cdot \bar g = 0$.
    Hence $\{a\}\cdot g=2h$ for some $h\in K_n F$ and $N_{E/F}
    \gamma=2(h+f)$. Thus
    \begin{equation*}
        N_{E/F}\left(\gamma-i_E(h+f)\right)=0.
    \end{equation*}
    Then by Theorem~\ref{th:bk} there exists $\alpha\in K_n E$
    such that
    \begin{equation*}
        (\sigma-1)\alpha=\gamma-i_E(h+f).
    \end{equation*}
    Observe that since $\bar\gamma\notin i_E k_n F$ we have
    $\bar\gamma-i_E(\overline{h+f})\neq 0$. Hence
    \begin{equation*}
        \bar\gamma = (\sigma-1)\bar\alpha + i_E(\overline{h+f})
        \in N k_n E,
    \end{equation*}
    as required.

    (2).  Suppose $\bar\gamma\in (k_nE)^G$.  Then $l(\gamma)=1$
    and the preceding part of our proof shows that $\langle \bar
    \gamma \rangle = \langle N\bar\alpha \rangle$ for $\alpha \in
    K_nE$.  Hence $(k_nE)^G\subset i_EN_{E/F}k_nE$.  Then
    \begin{equation*}
        i_EN_{E/F}k_nE \subset i_Ek_nF \subset (k_nE)^G \subset
        i_EN_{E/F}k_nE,
    \end{equation*}
    and so all inclusions are equalities.
\end{proof}

Our second lemma establishes a situation in which all elements in
$k_nE$ are fixed by $G$.

\begin{lemma}\label{le:extend2}
    Let $n\in \N$.  Suppose that
    \begin{equation*}
        i_E(\{\xi_p\}\cdot k_{n-1}F) = i_EN_{E/F}k_nE = \{0\}.
    \end{equation*}
    Then $(k_nE)^G=k_nE$.
\end{lemma}

\begin{remark*} The hypothesis $i_E(\{\xi_p\}\cdot k_{n-1}F)=
\{0\}$ can be omitted in the case $p=2$.
\end{remark*}

\begin{proof}
    Let $\gamma \in K_nE$.  We show that $l(\gamma)>1$ leads to
    a contradiction, whence we will have the result.

    Suppose that $l=l(\gamma)\ge 2$ and $1\le l\le i\le p$. We
    show by induction on $i$ that there exists $\alpha_i\in K_nE$
    such that $\langle (\sigma-1)^{i-1} \bar\alpha_i \rangle =
    \langle (\sigma-1)^{l-1} \bar\gamma \rangle$.  If $i=l$ then
    $\alpha_i=\gamma$ suffices.  Assume now that $l\le i<p$ and
    that our statement is true for $i$. Set $c=N_{E/F}\alpha_i$.
    Since $i_E\bar c = N\bar\alpha_i = (\sigma-1)^{p-1}
    \bar\alpha_i$ and $i<p$, $i_E\bar c = 0$.

    By Theorem~\ref{th:esext}, $\bar c = \{a\} \cdot \bar b$ for
    some $b\in K_nF$.  Hence $c = \{a\} \cdot b + pf$ for $f\in
    K_nF$.  Then since $2\le i<p$ in this case,
    \begin{equation*}
        N_{E/F}\left(\alpha_i-\{\root{p}\of{a}\}\cdot
        i_E(b)-i_E(f) \right) = 0.
    \end{equation*}
    By Theorem~\ref{th:bk}, there exists $\omega\in K_nE$ such
    that
    \begin{equation*}
        (\sigma-1)\omega = \alpha_i-\{\root{p}\of{a}\} \cdot
        i_E(b) - i_E(f).
    \end{equation*}
    Then $(\sigma-1)^2\omega = (\sigma-1)\alpha_i - i_E(\{\xi_p\}
    \cdot i_E(b)) = (\sigma-1)\alpha_i \neq 0$, and we can set
    $\alpha_{i+1}=\omega$. Observe that here we use our hypothesis
    \begin{equation*}
        i_E(\{\xi_p\}\cdot k_{n-1}F)=\{0\}.
    \end{equation*}
    Hence by induction there exists $\alpha_p\in K_nE$ such that
    \begin{equation*}
        \langle N\bar\alpha_p \rangle = \langle (\sigma-1)
        ^{l-1}\bar\gamma \rangle.
    \end{equation*}

But $i_EN_{E/F}\bar\alpha_p = 0$, whence $(\sigma-1)^{l-1} \bar\gamma
= 0$, a contradiction.
\end{proof}

Finally, we record a necessary and sufficient condition for an
$\Fp[G$]-module to be free.

\begin{lemma}\label{le:free}
    Let $M$ be an $\Fp[G]$-module.  Then the following are
    equivalent:
    \begin{enumerate}
        \item $M$ is a free $\Fp[G]$-module
        \item $M^G=(\sigma-1)^{p-1}M$.
    \end{enumerate}
\end{lemma}

\begin{proof}
    Condition (2) is equivalent to
    \begin{equation*}
        H^2(G,M) = \{0\}.
    \end{equation*}
    But this condition is known to be equivalent with (1) (for any
    $p$-group $G$!). (See for example \cite[p.~63]{La2}.)
\end{proof}

\subsection{Hereditary Properties}\

We say that a property of Milnor $k$-groups $k_nE$ and $k_nF$ is
\emph{hereditary} if the validity of the property for a given $n$
implies the validity of the property for all integers greater than
$n$.

The next lemma establishes various hereditary properties,
including the properties in Corollary~\ref{co:injsur}.

\begin{lemma}\label{le:normher}
    Let $n\in \N$.

    The following are hereditary properties:
    \begin{enumerate}
        \item $k_{n-1}F = \ann_{n-1} \{a\}  = \ann_{n-1}
        \{a,\xi_p\}.$
        \item $i_E\colon k_nF \to k_nE$ is injective
        \item $N_{E/F}\colon k_{n-1}E\to k_{n-1}F$ is surjective
        \item for some fixed $\alpha_1, \alpha_2 \in K_1F$,
        $\bar\alpha_1 \cdot k_{n-1}F \subset \bar\alpha_2 \cdot
        k_{n-1}F$
        \item for some fixed $\alpha\in K_1E$, $k_nE = i_Ek_nF +
        \bar\alpha \cdot i_Ek_{n-1}F$
    \end{enumerate}
\end{lemma}

\begin{proof}
    (1).  $k_nF = k_{n-1}F \cdot k_1F$, and since $\ann_{n-1}
    \{a\} = k_{n-1}F$, we have $\ann_n \{a\} = k_nF$ as well. The
    other equality follows from $\ann_{n} \{a\} \subset \ann_n
    \{a,\xi_p\}$.  The result follows by induction.

    (2-3). By Corollary~\ref{co:injsur}, the first three
    properties are equivalent, hence (2) and (3) are hereditary.

    (4). $K_nF = K_{n-1}F\cdot K_1F$, so $K_nF$ is generated by
    elements of the form
    \begin{equation*}
        \{f_1, f_2, \dots, f_n\} = \{f_1,\dots,f_{n-1}\} \cdot
        \{f_n\}, \quad f_i\in F^\times.
    \end{equation*}
    For each such generator, we calculate
    \begin{equation*}
        \bar\alpha_1 \cdot \overline{\{f_1,\dots,f_n\}} =
        \bar\alpha_2 \cdot \bar g \cdot \overline{\{f_n\}}
    \end{equation*}
    for some $g\in K_{n-1}F$, whence $\bar\alpha_1 \cdot k_nF
    \subset \bar\alpha_2 \cdot k_nF$.  The result follows by
    induction.

    (5). $k_{n+1}E = k_1E\cdot k_nE$, so the condition on $k_nE$
    gives us that $k_{n+1}E$ is generated by elements of the form
    \begin{equation*}
        \bar\gamma_1 = \overline{\{\delta\}} \cdot
        i_E(\overline{\{f_1,\dots,f_n\}}), \quad \delta\in
        E^\times, \ f_i\in F^\times
    \end{equation*}
    and
    \begin{equation*}
        \bar\gamma_2 = \overline{\{\delta\}} \cdot \bar\alpha
        \cdot i_E(\overline{\{f_1,\dots,f_{n-1}\}}), \quad
        \delta\in E^\times,\ f_i\in F^\times.
    \end{equation*}

    If $n-1\geq 1$ then we see that $k_{n+1}E$ is generated by the
    elements in $k_n E\cdot i_E k_1 F$. By hypothesis $k_nE = i_E
    k_nF + \bar\alpha \cdot i_E k_{n-1}F$ and therefore $k_{n+1}E$
    is generated by elements in $i_E k_{n+1}F + \bar\alpha\cdot
    i_E k_{n}F$.

    If $n=1$ then using our hypothesis $k_1 E = i_E k_1
    F + \bar\alpha \cdot i_E k_0 F$ we may write the generators
    $\bar\gamma_2$ of $k_2 E$ as
    \begin{equation*}
        \bar\gamma_2 = \left(i_E(\overline{\{f\}}) + c\bar\alpha
        \right)\cdot \bar \alpha, \quad f
        \in F^\times,\ c\in\Z.
    \end{equation*}
    Since $\bar\alpha\cdot\bar\alpha = \{-1\} \cdot \bar \alpha$,
    \begin{equation*}
        \bar\gamma_2 = i_E(\overline{\{f\}}) \cdot \bar\alpha +
        c\ i_E(\{-1\}) \cdot \bar\alpha = - \bar \alpha \cdot
        i_E(\overline{\{f\}}) - \bar\alpha \cdot i_E(c\{-1\}).
    \end{equation*}
    Thus in this case both types of generators of $k_2 E$ have the
    required form of elements in $i_E k_2 F+ \bar\alpha \cdot i_E
    k_1 F$.

    The result now follows by induction.
\end{proof}

\subsection{Fields of the Form {$\C((\oplus_I \Z_{(p)}))$}}\

For our examples in sections~\ref{se:exfree} and \ref{se:extriv}
we introduce the following notation and results.

Let
\begin{equation*}
    \Z_{(p)}:=\left\{\frac{c}{d}\in\Q\ \big\vert\  c,d\in\Z, d
    \neq 0; \ \text{if } c\neq 0 \text{ then } (c,d)=1, p\nmid
    d\right\}.
\end{equation*}
Observe that $\Z_{(p)}$ carries a natural ordering induced from
$\Q$.  Let $I$ be a well-ordered set of cardinality $m$, and let
$\Gamma$ be a direct sum of $m$ copies of $\Z_{(p)}$, indexed by
$I$. Then $m=\dim_{\Fp}\Gamma /p\Gamma$. Order $\Gamma$
lexicographically.

Then $\Gamma$ is a linearly ordered abelian group. (Recall that
each non-empty set can be well-ordered (see \cite[Appendix~2,
Theorem~4.1]{La}).) Now it is well-known that since $\Gamma$ is a
totally ordered abelian group, the field
\begin{equation*}
    F_m := \C((\Gamma)) :=\{f\colon \Gamma\to \C\ \vert\ \supp
    (f)\text{ is well-ordered}\}
\end{equation*}
is a henselian valued field with value group $\Gamma$ and residue
field $\C$. (See \cite[Chapitre D, Th\'eor\`emes 2 et 3, page 103,
et Chapitre F, Th\'eor\`eme 4, page 198]{Rib1}.) Thus a typical
element $f\in F_m$ may be written as a formal sum
\begin{equation*}
    f=\sum_{g\in\Gamma} a_g t^g
\end{equation*}
such that the set $\supp (f):=\{g\in \Gamma\ \vert\ a_g\neq 0\}$ is
a well-ordered subset of $\Gamma$. The absolute Galois group of
$F_m$ is known to be $\Z_p^m$, the topological product of $m$ copies
of $\Z_p$ \cite[pages~3~and~4]{K}.  We record one property of $F_m$
in the following lemma.

\begin{lemma}\label{le:fmlemma}
    For $m, n\in \N \cup \{\aleph_0\}$,
    \begin{equation*}
        H^n(F_m)\cong \bigwedge^n H^1(\Z_p^m) \cong \bigwedge^n
        \oplus_m \Fp
    \end{equation*}
    where the cup-product is sent to the wedge product.
\end{lemma}

\begin{proof}
    Since $F_m$ is a henselian valued field, the second result
    follows from \cite[Theorem 3.6]{W}, observing that under the
    Kummer isomorphism $F_m^\times/F_m^{\times p}\cong H^1(F_m)$,
    $0=(-1)\in H^1(F_m)$, and $H^j(\C)=\{0\}$ for all $j\in \N$.
\end{proof}

Of particular interest to us will be certain fields with absolute
Galois groups which are pro-$p$ free products of groups of the
form $\Z_p^m$.

\begin{lemma}\label{le:fmm}
    Suppose that $m_1, m_2$ are non-zero cardinal numbers, and let
    $F_{m_1}$ and $F_{m_2}$ be as above.  There exists a field
    $F_{m_1,m_2}$ of characteristic $0$, containing a primitive
    $p^2$th root of unity $\xi_{p^2}$, such that the absolute
    Galois group
    \begin{equation*}
        G_{F_{m_1,m_2}} \ \ \cong \ \ G_{F_{m_1}}
        \star_{\text{pro-}p} G_{F_{m_2}} \ \ \cong \ \
        \Z_p^{m_1} \star_{\text{pro-}p} \Z_p^{m_2},
    \end{equation*}
    where the free products are taken in the category of
    pro-$p$-groups, and the natural restriction maps
    \begin{equation*}
        \res_{\star} \colon H^n(F_{m_1,m_2})\to H^n(F_{m_1})\oplus
        H^n(F_{m_2})
    \end{equation*}
    are isomorphisms.
\end{lemma}

Note that we use the notation $\res_{\star}$ to distinguish this
restriction map from restriction maps $H^n(F)\to H^n(E)$.  For
$h\in H^n(F_{m_1,m_2})$, we will write $\res_\star h = h_1\oplus
h_2$.

\begin{proof}
    The existence of a field $F_{m_1,m_2}$ with $\chr(F_{m_1,m_2})
    = \chr(F_{m_1}) = \chr(F_{m_2}) = 0$ and the given absolute
    Galois group follows from \cite[Proposition~1.3]{EH}.

    Additionally using the construction of $F_{m_1,m_2}$ following
    \cite[proof of Proposition~1.3]{EH} we assume that $F_{m_1,m_2}$
    is the intersection of two henselian valued fields $(L_i,V_i),
    i=1,2$, with residue fields isomorphic to $F_{m_1}$ and
    $F_{m_2}$ respectively. Here $V_i$ is a henselian valuation on
    $L_i$. Then by Hensel's Lemma (see \cite[pages~12 and~13,
    condition (3)]{Rib2}) and by the fact that $F_{m_1}$ and
    $F_{m_2}$ have characteristic $0$ and both contain a primitive
    $p^2$th root of unity, we see that $F_{m_1,m_2}$ also contains a
    primitive $p^2$th root of unity. The fact that the restriction
    maps are isomorphisms follows from \cite[S\"atze (4.1) und
    (4.2)]{N}.
\end{proof}

\begin{remark*}
    From the proof above it follows that $F_{m_1,m_2}$ contains all
    $p^k$th primitive roots, $k\in\N$. However we shall not need
    this observation.
\end{remark*}

\section{When is Galois Cohomology Free?}\label{se:free}

\begin{proof}[Proof of Theorem~\ref{th:free}]
    Here and elsewhere we use Theorem~\ref{th:bk} to translate
    between Galois cohomology $H^n$ and $K$-theory $k_n$.

    First we show that for all $p$, $k_nE$ free implies that
    \begin{equation*}
        i_EN_{E/F}k_nE = i_Ek_nF = (k_nE)^G.
    \end{equation*}
    If $k_nE$ is free, then by Lemma~\ref{le:free},
    $(\sigma-1)^{p-1} k_nE = (k_nE)^G$. Observing that
    $(\sigma-1)^{p-1}$ and $i_EN_{E/F}$ are equivalent,
    \begin{equation*}
        i_EN_{E/F}k_nE = (\sigma-1)^{p-1}k_nE = (k_nE)^G.
    \end{equation*}
    Then since $i_Ek_nF\subset (k_nE)^G$ and $i_EN_{E/F}k_nE
    \subset i_Ek_nF$, we have established our claim.

    Assume first that $p>2$.  First we show (1)$\implies$(2). Let
    $f\in K_{n-1}F$ be arbitrary, and set $\alpha =
    \{\root{p}\of{a}\}\cdot f$.  Now because $(k_nE)^G =
    i_EN_{E/F}k_nE$, there exists $\beta\in K_nE$ such that
    $i_EN_{E/F}\bar\beta = i_E(\{\xi_p\}\cdot \bar f)\in
    (k_nE)^G$.  Set $\gamma = (\sigma-1)^{p-2}\beta$.  Since
    $p>2$, $\gamma$ is in the image of $\sigma-1$ and hence has
    trivial norm. We calculate
    \begin{equation*}
        N_{E/F}(\bar\alpha - \bar\gamma) = \{a\} \cdot \bar f.
    \end{equation*}
    On the other hand, observing that $i_EN_{E/F} =
    (\sigma-1)^{p-1}$ on $k_nE$,
    \begin{equation*}
        (\sigma-1)(\bar\alpha - \bar\gamma) = i_E(\{\xi_p\} \cdot
        \bar f) - i_E(\{\xi_p\} \cdot \bar f) = 0.
    \end{equation*}
    Hence $\bar\alpha - \bar\gamma \in (k_nE)^G = i_Ek_nF$.  But
    on $i_Ek_nF$ the norm map $N_{E/F}$ is trivial.  Hence $\{a\}
    \cdot \bar f = 0$ and $\ann_{n-1} \{a\} = k_{n-1}F$.

    By Corollary~\ref{co:injsur}, (2), (3), and (4) are all
    equivalent.  Now we show (4)$\implies$(1). Assume that
    $N_{E/F}\colon k_{n-1}E\to k_{n-1}F$ is surjective.  By
    Lemma~\ref{le:extend} we have $(k_nE)^G = (\sigma-1)^{p-1}
    k_nE$. Hence by Lemma~\ref{le:free}, $k_nE$ is free.

    Now suppose that $p=2$.  By Corollary~\ref{co:cop2}, we need
    only show that (1) and (2) are equivalent.  We show first that
    (1)$\implies$(2). We established that (1) implies
    $i_EN_{E/F}k_nE = i_Ek_nF$.  Since $\ker i_E = \{a\} \cdot
    k_{n-1}F$, this equality is equivalent to $k_nF = N_{E/F}k_nE
    + \{a\} \cdot k_{n-1}F$, so we have the second part of (2).
    Clearly $\ann_{n-1} \{a\} \subset \ann_{n-1} \{a,-1\}$, so we
    show that $\ann_{n-1} \{a,-1\} \subset \ann_{n-1} \{a\}$.

    We adapt the argument above.  Let $\bar f \in \ann_{n-1} \{a,
    -1\}$.  Set $\bar\alpha = \{\sqrt{a}\} \cdot \bar f$.  Since
    $\{a\}\cdot \{-1\}\cdot \bar f = 0$, Theorem~\ref{th:esext}
    tells us that there exists $\beta \in K_{n}E$ such that
    $N_{E/F}\bar\beta = \{-1\}\cdot \bar f$.  Now we calculate by
    the projection formula
    \begin{equation*}
        N_{E/F}(\bar\alpha - \bar\beta) = \{-a\} \cdot \bar f
        - \{-1\} \cdot \bar f = \{a\} \cdot \bar f.
    \end{equation*}
    On the other hand, using the fact that $\sigma-1=\sigma+1$
    when $p=2$,
    \begin{equation*}
        (\sigma-1)(\bar\alpha - \bar\beta) = \{-1\} \cdot \bar f -
        \{-1\} \cdot \bar f = 0.
    \end{equation*}
    Hence $\bar\alpha-\bar\beta\in (k_nE)^G$.  By
    Lemma~\ref{le:free}, $(k_nE)^G=(\sigma-1)k_nE =
    i_EN_{E/F}k_nE\subset i_Ek_nF$. Therefore $\bar\alpha -
    \bar\beta \in i_Ek_nF$.  But on $i_Ek_nF$, the norm map
    $N_{E/F}$ is trivial.  Hence $\{a\} \cdot \bar f = 0$, and
    $\bar f\in \ann_{n-1} \{a\}$, so $\ann_{n-1} \{a,-1\} \subset
    \ann_{n-1} \{a\}$, as required.

    Now we show that (2)$\implies$(1). Assume that $\ann_{n-1}
    \{a,-1\}$ $=$ $\ann_{n-1} \{a\}$ and that $k_nF = N_{E/F}k_nE +
    \{a\} \cdot k_{n-1}F$.  We use Lemmas~\ref{le:extend} and
    \ref{le:free} to deduce that $k_nE$ is free.
\end{proof}

It follows easily that
\begin{corollary}\label{co:h1notfree}
    For $p>2$, $k_1E$ is never free.
\end{corollary}

\begin{proof}
    Since $G$ acts trivially on $k_0E\cong \Fp$,
    \begin{equation*}
        N_{E/F}k_0E = 0 \neq k_0F \cong \Fp.
    \end{equation*}
    Alternatively, $i_E\colon k_1F\to k_1E$ is not injective,
    since $\{a\}\in k_1F$ is a nontrivial element of the kernel.
\end{proof}

With Theorem~\ref{th:free} in hand, Lemma~\ref{le:normher} is
enough to establish hereditary freeness in the $p>2$ case, and for
the $p=2$ case we show that an additional condition, analogous to
the $p>2$ case, is sufficient:

\begin{corollary}\label{co:freeherp2}
    Suppose that $p=2$ and for some $n\in \N$,
    \begin{equation*}
        \ann_{n-1} \{a\} = k_{n-1}F.
    \end{equation*}
    Then $k_mE$ is a free $\F_2[G]$-module for all $m\ge n$.
\end{corollary}

\begin{proof}
    We show that the two conditions of part (2) of the $p=2$
    portion of Theorem~\ref{th:free} hold for $K$-theory degree at
    least $n$.  From Lemma~\ref{le:normher}, part (1), we deduce
    that $k_mF = \ann_m \{a\} = \ann_m \{a,-1\}$ for all $m\ge
    n-1$.

    By Theorem~\ref{th:esext} and Lemma~\ref{le:normher}, we have
    $k_mF = N_{E/F}k_mE$ for all $m\ge n-1$ and therefore we see
    that
    \begin{equation*}
        k_mF = N_{E/F}k_mE + \{a\} \cdot k_{m-1}F \quad
        \text{for all } m\ge n.
    \end{equation*}

    We conclude that $k_mE$ is a free $\F_2[G]$-module for all
    $m\ge n$.
\end{proof}

Just as before it follows easily that
\begin{corollary}\label{co:h1notfreep2}
    For $p=2$ and $\sqrt{-1}\in F$, $k_1E$ is never free.
\end{corollary}

\begin{proof}
    Since $-1\in F^{\times 2}$, we have $\{-1\} = 0\in k_1F$ and
    $\{a,-1\}=0\in k_2F$, so that $\ann_0 \{a,-1\} = k_0F \cong
    \F_2 \neq \ann_0 \{a\} = \{0\}$.
\end{proof}

We are now ready to prove Theorem~\ref{th:freeher}.

\begin{proof}[Proof of Theorem~\ref{th:freeher}]
    For $p>2$, the fact that free cohomology is hereditary follows
    from Lemma~\ref{le:normher} and condition (2) in
    Theorem~\ref{th:free}.  The exactness of the first term of the
    sequence follows from Theorem~\ref{th:free}, part (3), while
    the exactness at the third term follows from
    Theorem~\ref{th:free}, part (4) and Lemma~\ref{le:normher},
    part (3).

    Assume then that $p=2$ and $a=x^2+y^2$ for some $x$, $y\in
    F^\times$. If $-1\in F^{\times 2}$ then $\{a,-1\}\in 2K_2F$
    and so $\{a,-1\} = 0\in k_2F$.  Otherwise let
    $K=F(\sqrt{-1})$, and observe that $a=N_{K/F}(x+y\sqrt{-1})$.
    Then $\{a,-1\} = 0 \in k_2F$.  Hence $\ann_{n-1} \{a,-1\} =
    k_{n-1}F$.

    Now observe that since $k_nE$ is a free $\F_2[G]$-module, by
    Theorem~\ref{th:free} we have $\ann_{n-1} \{a\} = \ann_{n-1}
    \{a,-1\}$, and so $\ann_{n-1} \{a\} = \ann_{n-1} \{a,-1\} =
    k_{n-1}F$. We deduce from Corollary~\ref{co:freeherp2} that
    $k_{m}E$ is a free $\F_2[G]$-module for all $m\ge n$.

    For the exact sequence in the case $p=2$, we have shown that
    $k_{m-1}F = \ann_{m-1} \{a\}$, and so by
    Theorem~\ref{th:esext} and Lemma~\ref{le:normher}, we have
    $k_nF = N_{E/F}k_nE$ for all $n\ge m-1$.  Hence we have
    exactness at the third term.  Furthermore, we conclude from
    Corollary~\ref{co:injsur} that $i_E$ is injective from
    $k_{m}F$ to $k_{m}E$. Hence we have exactness at the first
    term as well.
\end{proof}

We now provide an example of $k_1E$ free but $k_2E$ nonfree,
showing that freeness is not generally hereditary when $p=2$.

\begin{example*}
    Let $p=2$, $F=\Q_2$, and $a=-1$, so $E=\Q_2(\sqrt{-1})$.
    Then
    \begin{equation*}
        k_1 F\cong F^\times/F^{\times 2}= \langle
        [-1],[2],[5]\rangle.
    \end{equation*}
    and
    \begin{equation*}
        N_{E/F}k_1 E \cong N_{E/F}(E^\times)
        F^{\times 2}/F^{\times 2} = \langle [2],[5]
        \rangle.
    \end{equation*}
    (See \cite[page~162, Corollary~2.24]{Lam}.)

    Therefore $k_1 F=N_{E/F}k_1 E+\{-1\}\cdot k_0 F$.
    Moreover, since $[-1]\notin N_{E/F}\left(E^\times\right)
    F^{\times 2}/F^{\times 2}$, we have $\{-1,-1\}\neq 0\in k_2
    F$. (Again see \cite[page~162, Corollary~2.24]{Lam}.)

    Hence $\ann_0 \{-1,-1\}=\{0\}$. Since $\{-1\}\neq 0$ in $k_1
    F$ we see that $\ann_0 \{-1\}=\{0\}$. Hence the conditions of
    part (2) of the $p=2$ portion of Theorem~\ref{th:free} are
    satisfied, whence $k_1 E$ is a free $\F_2$-module.

    Observe, however, that since $\ann_0 \{-1,-1\}=\{0\} \neq k_0
    F$, the first hypothesis in Corollary~\ref{co:freeherp2} does
    not hold. Therefore we cannot conclude that $k_2 E$ is a free
    $\F_2[G]$-module---and of course it is not, as it is well
    known that $k_2 E\cong\F_2$. (See \cite[page~158,
    Corollary~2.15]{Lam} and \cite[Theorem~4.1]{Mi}.)
\end{example*}

\section{Examples of $H^k(E)$ Free for all $n<k$ and $H^n(E)$
Nonfree, with Given Cohomological Dimension}\label{se:exfree}

We have shown in Theorem~\ref{th:freeher} that if $p>2$ then the
property $H^n(E)$ is a free $\Fp[G]$-module is hereditary.
Moreover, the same property is hereditary in the case $p=2$ as
well if $i=\sqrt{-1}\in F$, since then $a = ((a+1)/2)^2 +
((a-1)i/2)^2$.

These results lead naturally to the definition of an interesting
invariant $\cf(E/F) \in \{0\} \cup \N \cup\{\infty\}$:
\begin{multline*}
    \cf(E/F) \ = \\ \sup \left\{ n\in\N\cup\{0\} \ \ \vert \ \
    H^n(E) \text{ is not a free } \Fp[G]\text{-module}\right\}.
\end{multline*}
We have chosen $\cf$ to indicate that after degree $\cf(E/F)$,  Galois
cohomology is cohomologically free.  Of course, if $H^n(E)$ is
never free then $\cf(E/F)=\infty$, and otherwise $\cf(E/F) \in
\N\cup \{0\}$.

Assume for the moment that either $p>2$ or $\sqrt{-1}\in F$. If
$\cf(E/F)=n$, then by definition $H^m(E)$ is a free
$\Fp[G]$-module for all $m>n$. On the other hand, by the
hereditary property we also have that $H^k(E)$ is not free for all
$k\leq \cf(E/F)$. Finally, Corollaries~\ref{co:h1notfree} and
\ref{co:h1notfreep2} tell us that $H^1(E)$ is never free and hence
$\cf(E/F)\geq 1$. A natural question arises: can we choose a
suitable field extension $E/F$ so that $\cf(E/F)$ is a given
natural number or $\infty$? We show that the answer is
affirmative.

Before formulating our result precisely, let us recall that for
any pro-$p$-group $T$ we may define $\cd(T)$, the cohomological
dimension of $T$, as
\begin{equation*}
    \cd(T) \ = \ \sup \left\{ k \in \N \ \ \vert \ \ H^k(T,\Fp)
    \neq \{0\} \right\} \ \in \ \N \cup \{\infty\}.
\end{equation*}
(See \cite[Chapter 7]{RZ}.) Suppose that the absolute Galois group
$G_E$ is a pro-$p$-group. Then, adopting the convention that
$\{0\}$ is considered a free $\Fp[G]$-module,  we have:
\begin{equation*}
    \cf(E/F) \leq \cd(G_E).
\end{equation*}

From Corollaries~\ref{co:h1notfree} and \ref{co:h1notfreep2} above
it follows that if $p>2$ or $p=2$ and $\sqrt{-1}\in F$ then
$\cf(E/F)\ge 1$.

Our result is then the following.

\textit{Given $1\le n\le m \in \N \cup \{\infty\}$ and a prime
$p$, there exists a cyclic extension $E/F$ of degree $p$ with
$\xi_p \in F$ such that
\begin{enumerate}
    \item $G_E$ is a pro-$p$-group;
    \item $\cf(E/F)=n$; and
    \item $\cd(G_E)=m$.
\end{enumerate}
}

Observe that if we choose $n\in\N$ and $m>n$, then we have
obtained examples as promised in the title of this section.

\subsection{The case $m\in \N$.}\

(1). Let $F:=F_{n,m}$ be a field of characteristic $0$ with
$G_{F}\cong \Z_p^n \star_{\text{pro-}p} \Z_p^{m}$ and $\xi_{p^2}\in
F$, given by Lemma~\ref{le:fmm}. Observe particularly that
$\sqrt{-1}\in F$ in the case $p=2$.  Let
\begin{equation*}
    E=F(\root{p}\of{a})
\end{equation*}
for any $a\in F^\times$ such that under the restriction map on
$H^1$,
\begin{equation*}
    \res_{\star} (a) = (a)_1\oplus (a)_2, \quad (a)_1\neq 0, \
    (a)_2 = 0.
\end{equation*}
We use here, and later without mention, the fact that
$\res_{\star}$ is an isomorphism, by Lemma~\ref{le:fmm}. Observe
that there exists an $a$ with the required conditions because by
Lemma~\ref{le:fmlemma}, $H^1(F_n)\neq\{0\}$.

(2a). \textit{$H^n(E)$ is not free}. We claim that
\begin{equation*}
    \ann_{n-1} (a) \neq H^{n-1}(F).
\end{equation*}
If $n=1$ this inequality is true as $(a)\ne 0\in H^1(F)$. Assume
now that $n>1$. We shall use Lemma~\ref{le:fmlemma} together with
the fact that the restriction map in the cohomology ring of a
profinite group to the cohomology ring of a closed subgroup is
a ring homomorphism. (See for example \cite[Proposition 7.9.4]{RZ}.)

Let $a_1\in F_{n}^\times$ satisfy $(a_1)=(a)_1$, and extend
$\{(a_1)\}$ to a basis $\{(a_1)$, $(a_2)$, $\cdots$, $(a_n)\}$ of
$H^1(F_n)$.  By Lemma~\ref{le:fmlemma}, the element
\begin{equation*}
    (a_1)\cup (a_2) \cup \cdots \cup (a_n) \in H^n(F_n)
\end{equation*}
is nontrivial, so that $0\neq (a_2)\cup \cdots \cup (a_n)\in
H^{n-1}(F_n)$. Let $b\in H^{n-1}(F)$ satisfy
\begin{equation*}
    b_1= (a_2)\cup \cdots \cup (a_n)\in H^{n-1}(F_n), \qquad
    b_2=0\in H^{n-1}(F_{m}).
\end{equation*}
Then since the cup-product commutes with $\res_{\star}$,
\begin{equation*}
    \big((a)\cup b\big)_1 = (a_1)\cup b_1 \neq 0\in H^n(F_n),
\end{equation*}
so that $(a)\cup b\neq 0\in H^n(F)$ and hence $\ann_{n-1} (a) \neq
H^{n-1}(F)$.

If $p>2$, we conclude by Theorem~\ref{th:free} that $H^n(F)$ is
not free.  If $p=2$, observe that since $\sqrt{-1}\in F$, we have
$\ann_{n-1} (a,-1) = \ann_{n-1} 0 = H^{n-1}(F)$, so that
$\ann_{n-1} (a,-1) \neq \ann_{n-1} (a)$.  We deduce from
Theorem~\ref{th:free} that $H^n(F)$ is not free.

(2b). \textit{$H^{k}(E)$ is free for all $k\ge n+1$}. We claim
that
\begin{equation*}
    \ann_n (a) = H^n(F).
\end{equation*}

Let $c\in H^n(F)$.  Then since $H^{n+1}(F_n)=0$ by
Lemma~\ref{le:fmlemma},
\begin{equation*}
    \res_{\star} (a)\cup c = \big((a_1) \cup c_1\big) \oplus \big(0
    \cup c_2\big) = 0 \oplus 0 = \res_{\star} 0
\end{equation*}
Hence $(a)\cup c=0$ and $\ann_{n} (a) = H^{n}(F)$.

If $p>2$ then we conclude by Theorem~\ref{th:free} that
$H^{n+1}(E)$ is free, and by Theorem~\ref{th:freeher}, $H^k(E)$ is
free for all $k\ge n+1$. If $p=2$, observe that $\ann_{n} (a,-1) =
\ann_n 0 = H^n(F)$.  Furthermore, we use Corollary~\ref{co:injsur}
to obtain that $\cor\colon H^{n}(E)\to H^{n}(F)$ is surjective.
Then by Corollary~\ref{co:freeherp2}, we have that $H^{k}(E)$ is
free for all $k\ge n+1$.

(3). \textit{$\cd(G_E)=m$}. First we claim that $G_E$ does not
contain an element of order $p$. By Artin-Schreier's theorem (see
for instance \cite[Chapter VI, Theorem 17]{J}), finite subgroups
of absolute Galois groups are either trivial or of order $2$, and
since $\sqrt{-1}\in E$ no element of order $2$ exists in $G_E$.

Then, by Serre's well-known theorem \cite{S}, we obtain
\begin{equation*}
    \cd(G_E)  = \cd(G_F).
\end{equation*}
From Lemmas~\ref{le:fmlemma} and \ref{le:fmm} we find that
\begin{equation*}
    \cd(G_F) = \max \{ \cd(F_n), \cd(F_m)\} = m.
\end{equation*}
Thus $\cd(G_E)=m$ as required.

\subsection{The case $n<m=\infty$.}\label{ss:nminfty}
Set
\begin{equation*}
    F_\infty := \C\left( \left(\Z_{(p)}^m\right) \right), \quad
    \text{ where }m = \aleph_0.
\end{equation*}
With the same argument as in the proof of Lemma~\ref{le:fmm},
there exists a field $F:=F_{n,\infty}$ such that $G_F\cong
G_{F_{n}} \star_{ \text{pro-}p} G_{F_{\infty}}$ and $\xi_{p^2} \in
F$. Then set
\begin{equation*}
    E=F(\root{p}\of{a})
\end{equation*}
for any $a\in F^\times$ such that under the restriction map
\begin{equation*}
    \res_\star \colon H^1(F)\to H^1\left(F_{n}
    \right)\oplus H^1\left(F_\infty\right),
\end{equation*}
we have
\begin{equation*}
    \res_\star (a)=(a)_1\oplus 0,\quad (a)_1\neq 0.
\end{equation*}

Then $\cd(G_F)=\cd(G_E)=\infty$, and with the same argument as
above we see that $\cf(E/F)=n$.

\subsection{The case $n=\infty=m$.}  As above we let $\Gamma$
be a direct sum of $\aleph_0$ copies of $\Z_{(p)}$. Then we set
$F:=F_\infty = \C((\Gamma))$. Let $a\in F^\times$ such that
$v(a) \in \Gamma\setminus p\Gamma$, where $v$ is a natural
valuation on $F$. Then from the description of Galois cohomology
of $p$-henselian fields (see \cite[Theorem~3.6]{W}), we obtain
\begin{equation*}
    \ann_n (a) = (a) \cup H^{n-1}(F)
\end{equation*}
and
\begin{equation*}
    (a) \cup H^{n-1}(F) \neq H^n(F) \quad \text{for all } n \in
    \N.
\end{equation*}
(Observe that when $p=2$ we use the fact that $\sqrt{-1}\in F$ in
the cited result.) Setting $E=F(\root{p}\of{a})$, just as before
we have that
\begin{equation*}
    \cf(E/F)=\infty,
\end{equation*}
as required.

\section{When is Galois Cohomology Trivial?}\label{se:trivial}

First we need a lemma.

\begin{lemma}\label{le:spectrivp2}
    Suppose that $p=2$. Then
    \begin{equation*}
        \{a\} \cdot \ann_{n-1} \{a,-1\} \subset N_{E/F}
        k_n E.
    \end{equation*}
\end{lemma}

\begin{proof}
    Let $\bar\beta \in \ann_{n-1} \{a,-1\}$. Then $\{-1\} \cdot
    \bar\beta \in \ann_n \{a\} = N_{E/F}k_n E$ by
    Theorem~\ref{th:esext}. Let $\gamma\in K_n E$ such that
    $\{-1\}\cdot \bar\beta = N_{E/F}(\bar\gamma)$. Then we have
    \begin{equation*}
        \{a\} \cdot \bar\beta = N_{E/F}(\{\sqrt{a}\}\cdot
        i_E(\bar\beta) + \bar\gamma).
    \end{equation*}
    Thus $\{a\}\cdot\ann_{n-1} \{a,-1\}\subset N_{E/F} k_n E$ as
    asserted.
\end{proof}

It is worth observing that if $n=1$, Lemma~\ref{le:spectrivp2} is
equivalent to
    \begin{equation*}
        \{a,-1\}=0 \quad \mbox{ if and only if } \quad \{a\}\in
        N_{E/F} k_1 E,
    \end{equation*}
and therefore Lemma~\ref{le:spectrivp2} can be viewed as a
generalization of this statement.

Now we are ready to prove Theorem~\ref{th:trivial}.

\begin{proof}[Proof of Theorem~\ref{th:trivial}]
    As before, we translate to $K$-theory using
    Theorem~\ref{th:bk}. We first consider the case $p>2$.

    (1)$\implies$(3). Assume that $k_n E$ is a trivial
    $\Fp[G]$-module. Suppose $f\in K_{n-1}F$, and set
    $\beta=\{\root{p}\of{a}\} \cdot i_E(f)$. Then
    \begin{equation*}
        (\sigma-1)\bar\beta=0 \implies  \{\xi_p\}\cdot i_E(\bar f)
        = i_E(\{\xi_p\}\cdot \bar f) = 0.
    \end{equation*}
    But then by Theorem~\ref{th:esext}, $\{\xi_p\}\cdot \bar f\in
    \{a\}\cdot k_{n-1}F$.

    Now let $\gamma\in K_nE$ be arbitrary.  Again, $(\sigma-1)
    \bar\gamma = 0$. Then
    \begin{equation*}
        i_EN_{E/F} \bar\gamma = (\sigma-1)^{p-1} \bar\gamma = 0
    \end{equation*}
    and so by Theorem~\ref{th:esext}, $N_{E/F}\bar\gamma =
    \{a\}\cdot \bar f$ for $f\in K_{n-1}F$.  By the projection
    formula, $N_{E/F}(\{\root{p}\of{a}\}\cdot i_E(\bar f)) =
    \{a\}\cdot \bar f$.  Then
    \begin{equation*}
        N_{E/F}(\bar\gamma - \{\root{p}\of{a}\} \cdot i_E(\bar
        f)) = 0,
    \end{equation*}
    and hence
    \begin{equation*}
        N_{E/F}(\gamma - \{\root{p}\of{a}\}\cdot i_E(f)) = p g,
    \end{equation*}
    for some $g\in K_{n-1}F$.  Set
    \begin{equation*}
        \beta = \gamma - \{\root{p}\of{a}\} \cdot i_E(f) -
        i_E(g).
    \end{equation*}
    Then $N_{E/F}(\beta)=0$.  By Theorem~\ref{th:bk}, there exists
    $\alpha\in K_nE$ such that $(\sigma-1)\alpha=\beta$.  But
    since $k_nE$ is fixed by $G$, $\bar\beta=0$.  Hence $k_nE =
    i_E(k_nF) + \{\root{p}\of{a}\}\cdot i_E(k_{n-1}F)$.

    (3)$\implies$(2).  Since $p>2$, $N_{E/F}(\{\root{p}\of{a}\}
    \cdot i_E(\bar f)) = \{a\} \cdot \bar f$ for $f\in K_{n-1}F$,
    and $N_{E/F}(i_E(\bar g))=0$ for $g\in K_nF$.  Hence $N_{E/F}
    k_nE=\{a\}\cdot k_{n-1}F$.  Since by Theorem~\ref{th:esext},
    $\ann_n \{a\} = N_{E/F} k_nE$, we are done.

    (2)$\implies$(1).  Assume that $\{\xi_p\} \cdot k_{n-1}F
    \subset \{a\}\cdot k_{n-1}F$ and $\{a\} \cdot k_{n-1}F =
    \ann_n \{a\}$.  By Theorem~\ref{th:esext}, $\ann_n \{a\} =
    N_{E/F}k_nE$.  Hence $\{\xi_p\} \cdot k_{n-1}F \subset
    N_{E/F}k_nE = \{a\} \cdot k_{n-1}F$.  But by
    Theorem~\ref{th:esext}, $\{a\} \cdot k_{n-1}F = \ker i_E$.  We
    then apply Lemma~\ref{le:extend2} to deduce that
    $k_nE=(k_nE)^G$.

    Now we consider the case $p=2$.

    (1)$\implies$(2). Assume that $k_n E$ is a trivial
    $\F_2[G]$-module. Let $\alpha\in K_nE$.  Then $i_EN_{E/F}
    \bar\alpha = (\sigma-1) \bar\alpha = 0$ implies that
    $N_{E/F}\bar\alpha = \{a\}\cdot \bar b$ for some $b\in
    K_{n-1}F$, by Theorem~\ref{th:esext}.  Now $\{a,-1\} =
    \{a,a\}$ in $k_2F$, and then
    \begin{equation*}
        \{a,-1\}\cdot \bar b = \{a\}\cdot N_{E/F} \bar\alpha = 0,
    \end{equation*}
    again by Theorem~\ref{th:esext}. Hence $N_{E/F}k_nE \subset
    \{a\} \cdot \ann_{n-1} \{a,-1\}$.  By Theorem~\ref{th:esext},
    $N_{E/F}k_nE = \ann_n \{a\}$, whence $\ann_n \{a\} \subset
    \{a\} \cdot \ann_{n-1} \{a,-1\}$.

    (2)$\implies$(1). Assume $\ann_n \{a\} \subset \{a\}\cdot
    \ann_{n-1} \{a,-1\}$.  By Theorem~\ref{th:esext}, we have
    $\ann_n \{a\} = N_{E/F}k_nE$ and therefore $N_{E/F}k_nE
    \subset \{a\} \cdot \ann_{n-1} \{a,-1\}$.   By
    Lemma~\ref{le:spectrivp2}, $\{a\} \cdot \ann_{n-1} \{a,-1\}
    \subset N_{E/F}k_nE$ and hence $N_{E/F} k_nE = \{a\} \cdot
    \ann_{n-1} \{a,-1\}$.  Let $\gamma\in K_nE$ be arbitrary. Then
    $N_{E/F}\bar\gamma = \{a\}\cdot \bar b$ for some $\bar b\in
    \ann_{n-1} \{a,-1\}$. Hence
    \begin{equation*}
        (\sigma-1) \bar\gamma = i_EN_{E/F} \bar\gamma =
        i_E(\{a\}\cdot \bar b).
    \end{equation*}
    But by Theorem~\ref{th:esext}, $i_E(\{a\}\cdot \bar b) = 0$.
    Hence $(\sigma-1) \bar\gamma = 0$, and $(k_nE)^G = k_nE$ as
    required.

    Now assume $p=2$ and $a\in (F^{\times 2} + F^{\times 2})
    \setminus F^2$. As in the proof of Theorem~\ref{th:freeher},
    we have that $\{a,-1\}=0\in k_2F$.  Therefore $\{-1\} \in
    \ann_1 \{a\}$. Since $\ann_1 \{a\} = N_{E/F}k_1E$ by
    Theorem~\ref{th:esext}, we obtain $\{-1\} \in N_{E/F}k_1E$.
    Equivalently, $-1\in N_{E/F}(E^\times)$. By \cite[Theorem
    3]{A}, $E/F$ embeds in an extension $E'/F$ cyclic of degree 4
    with $E'=E(\sqrt{\delta})$ for $\delta\in E^\times$. Kummer
    theory tells us that $\overline{\{\delta\}}\in (k_1E)^G$, so
    $(\sigma-1) (\overline{\{\delta\}})= (\sigma+1)
    (\overline{\{\delta\}})=0$. Therefore $(\sigma+1)(\{\delta\})
    \in 2K_1E$, whence $N_{E/F}(\delta)\in F^\times \cap E^{\times
    2}$.  On the other hand, Kummer theory gives that
    $(F^\times \cap E^{\times 2})/F^{\times 2} =  \{F^{\times 2},
    aF^{\times 2}\}$. If $N_{E/F}(\delta)=f^2$ for $f\in F^\times$
    then we have:
    \begin{align*}
    & \left(\sqrt{\delta}\right)^{\sigma^2-1}=\left(\sqrt{\delta}
    \sqrt{\delta}^{\sigma}\right)^{\sigma-1} = \left(\sqrt{\delta}
    \sqrt{\delta^\sigma}\right)^{\sigma-1} = \left(\sqrt{N_{E/F}
    \delta}\right)^{\sigma-1} \\ & = \left(\pm f\right)^{\sigma-1} = 1.
    \end{align*}
    (The choice of sign in the square roots above is irrelevant as
    we apply $\sigma-1$ afterwards.) Hence $\sigma$  extends to an
    order $2$ automorphism of $E'/F$, a contradiction.  We may
    conclude that $\overline{\{\delta\}} \in (k_1E)^G$ satisfies
    $N_{E/F} \overline{\{\delta\}} = \{a\}$, as required.

    We now follow the proof of the $p>2$ case to show that (1) and
    (3) are equivalent in the $p=2$ case as well.

    (1)$\implies$(3).  Assume that $k_n E$ is a trivial $\F_2[G]$-module.
    Let $\gamma\in K_nE$ be arbitrary.  Then $(\sigma-1) \bar\gamma =
    0$. Hence
    \begin{equation*}
        i_EN_{E/F} \bar\gamma = (\sigma-1) \bar\gamma = 0
    \end{equation*}
    and so by Theorem~\ref{th:esext}, $N_{E/F}\bar\gamma =
    \{a\}\cdot \bar f$ for $f\in K_{n-1}F$.  By the projection
    formula, $N_{E/F}(\{\delta\}\cdot i_E(\bar f)) = \{a\}\cdot
    \bar f$.  Then
    \begin{equation*}
        N_{E/F}(\bar\gamma - \overline{\{\delta\}} \cdot i_E(\bar
        f)) = 0,
    \end{equation*}
    and hence
    \begin{equation*}
        N_{E/F}(\gamma - \{\delta\} \cdot i_E(f)) = 2 g,
    \end{equation*}
    for some $g\in K_{n-1}F$.  Set
    \begin{equation*}
        \beta = \gamma - \{\delta\} \cdot i_E(f) -
        i_E(g).
    \end{equation*}
    Then $N_{E/F}(\beta)=0$.  By Theorem~\ref{th:bk}, there exists
    $\alpha\in K_nE$ such that $(\sigma-1)\alpha=\beta$.  But
    since $k_nE$ is fixed by $G$, $\bar\beta=0$.  Hence $k_nE =
    i_E(k_nF) + \overline{\{\delta\}}\cdot i_E(k_{n-1}F)$.

    (3)$\implies$(1). If $\beta\in K_nE$, then $\bar \beta =
    i_E(\bar g) + \overline{\{\delta\}} \cdot i_E(\bar f)$ for
    $g\in K_{n}F$ and $f\in K_{n-1}F$.  Clearly
    $(\sigma-1)\bar\beta = 0$, which implies that $(k_nE)^G =
    k_nE$.
\end{proof}

It is worth looking more closely at the case $n=1$ of
Theorem~\ref{th:trivial}.

Observe that the condition
    \begin{equation*}
        (\xi_p)\cup H^0(F)\subset (a)\cup H^0(F),
    \end{equation*}
for $n=1$ and $p>2$, is equivalent with the condition $\xi_
{p^2}\in E^\times$. (See \cite[Corollary~1]{MS}.)

In the case $n=1$ and $p=2$ the condition
    \begin{equation*}
        \ann_1(a)\subset (a)\cup \ann _0 (a,-1)
    \end{equation*}
can be reformulated as follows:
    \begin{equation*}
        \mbox{ If }(a,-1)=0\mbox{ then }N_{E/F}
        (k_1 E)\subset \left<\{a\}\right>
    \end{equation*}
and if $(a,-1)\neq 0$ then $\ann_1(a)=\{0\}$.

Since $\{-a,a\}=0$ we see that $\ann_1(a)=0$ implies that
$\{-a\}=0$ or equivalently $\{a\}=\{-1\}$. Recall that a
field $F$ is called Pythagorean if $F^2+F^2\subset F^2$.
In the first case when $(a,-1)=0$ the equality $\{-a,a\}
=0$ and $N_{E/F}(k_1 E)\subset \left<\{a\}\right>$ implies
$\sqrt{-1}\in E^\times$.

Summarizing our discussion for $p=2$ we have:

$H^1(E)$ is a trivial $\F_2[G]$-module if and only if either
$(a,-1)=0\in H^2(F)$ and $N_{E/F}(k_1 E)\subset \left<\{a\}
\right>$ or $a=-1$ and $F$ is a Pythagorean field. In both
cases $\sqrt{-1}\in E^\times$. (See \cite[Corollary~1]{MS}.)

\begin{corollary}\label{co:trives}
    Suppose $n\in \N$ and $(k_nE)^G=k_nE$. Then we have the
    following exact sequence:
    \begin{multline*}
        0 \to \ann_{n-1} \{a\} \to k_{n-1}F
        \xrightarrow{\{a\}\cdot -\ } k_nF \xrightarrow{i_E} \\
        k_nE \xrightarrow{N_{E/F}} \{a\} \cdot \ann_{n-1}
        \{a,\xi_p\} \to 0.
    \end{multline*}
    Here the map $\ann_{n-1} \{a\} \to k_{n-1}F$ is the natural
    inclusion.
\end{corollary}

\begin{proof}
    Exactness at the first and second terms is obvious, and
    exactness at the third term follows from
    Theorem~\ref{th:esext}.

    We consider exactness at the fifth term.  In the $p=2$ case,
    Theorem~\ref{th:trivial} tells us that $\ann_n \{a\} \subset
    \{a\} \cdot \ann_{n-1} \{a,-1\}$.  By Theorem~\ref{th:esext},
    we have $\ann_n \{a\} = N_{E/F}k_nE$, hence $N_{E/F}k_nE
    \subset \{a\} \cdot \ann_{n-1} \{a,-1\}$. By
    Lemma~\ref{le:spectrivp2} we have the reverse inclusion, so
    that $N_{E/F} k_nE = \{a\} \cdot \ann_{n-1} \{a,-1\}$ and the
    sequence is exact at the fifth term.

    In the $p>2$ case, observe that $\{\xi_p\}\cdot k_{n-1}F
    \subset \{a\} \cdot k_{n-1}F$ implies that $k_{n-1}F =
    \ann_{n-1} \{a,\xi_p\}$, since $\{a,a\}=0$. Therefore, by part
    (2) of Theorem~\ref{th:trivial}, we know $\ann_n \{a\} = \{a\}
    \cdot \ann_{n-1} \{a,\xi_p\}$. By Theorem~\ref{th:esext}, we
    have $\ann_n \{a\} = N_{E/F}k_nE$ and hence the sequence is
    exact at the fifth term in the $p>2$ case as well.

    Hence it remains to show exactness at the fourth term. Suppose
    $\gamma\in K_nE$ and $N_{E/F}\bar\gamma=0$.  Then there exists
    $f\in K_nF$ such that $N_{E/F}\gamma = pf$, and then
    $N_{E/F}(\gamma - i_E(f)) = 0$.  By Theorem~\ref{th:bk}, there
    exists $\alpha\in K_nE$ such that $(\sigma-1)\bar\alpha =
    \bar\gamma -i_E(\bar f)$.  But $(\sigma-1)\bar\alpha = 0$
    because $(k_nE)^G = k_nE$.  Hence $\bar\gamma = i_E(\bar f)$
    and we are done.
\end{proof}

We are now ready to prove Theorem~\ref{th:trivher}.

\begin{proof}[Proof of Theorem~\ref{th:trivher}]
    In the $p>2$ case, the result on heredity follows from
    Theorem~\ref{th:trivial}, part (3), together with two
    hereditary properties from Lemma~\ref{le:normher}: item (4),
    with $\alpha_1 = \{\xi_p\}$ and $\alpha_2 = \{a\}$, and
    item (5). The exact sequence, in turn, follows from
    Corollary~\ref{co:trives}.

    In the case $p=2$, by Theorem~\ref{th:trivial}, it is
    sufficient to prove that condition $(2)$ in the $p=2$ case is
    also hereditary. Assume $(2)$ holds for $n$ and $m>n$. By a
    well-known fact in Milnor $K$-theory, the group $K_m E$ is
    generated by the symbols
    \begin{align*}
        & \alpha = \{u,f_1,\dots,f_{n-1},\dots, f_{m-1}\}, \\
    & \mbox{ if }n>1\mbox{ and by }\alpha=\{u,f_1,\dots,
    f_{m-1}\}\mbox{ if }n=1,
    \end{align*}
    where $u\in E^*$ and $f_i\in F^*$ for all $i=1,\dots,m-1$.
    (See \cite[page~291, Corollary~2]{FV}.)

    Assume now that $n>1$. By the projection formula, we obtain
    \begin{equation*}
        N_{E/F}\bar\alpha = \overline{\{N_{E/F}u, f_1, \dots,
        f_{n-1}\}} \cdot \overline{\{f_n, \dots, f_{m-1}\}}.
    \end{equation*}
    Since $\ann_n \{a\} = N_{E/F}k_nE$ by Theorem~\ref{th:esext},
    condition (2) gives us that
    \begin{equation*}
        \overline{\{N_{E/F}u, f_1, \dots,
        f_{n-1}\}} = N_{E/F} \overline{\{u, f_1, \dots, f_{n-1}\}}
        \in \{a\}\cdot \ann_{n-1}\{a,-1\}.
    \end{equation*}
    Hence we may write
    \begin{equation*}
        \overline{\{N_{E/F}u, f_1, \dots,
        f_{n-1}\}} = \{a\}\cdot \bar c,
    \end{equation*}
    where $\bar c\in \ann_{n-1}\{a,-1\}$. Observe that this last
    equality holds also in the case when $n=1$, provided that we
    interpret the left-hand side as $\overline{\{N_{E/F}u\}}$.
    Thus
    \begin{equation*}
        N_{E/F}\bar\alpha = \{a\}\cdot\bar c\cdot \overline{\{ f_n,
        \dots, f_{m-1}\}}
    \end{equation*}
    and
    \begin{equation*}
        \bar c\cdot \overline{\{f_n,\dots, f_{m-1}\}} \in
        \ann_{m-1} \{a, -1\}.
    \end{equation*}
    Therefore $N_{E/F}k_m E \subset \{a\}\cdot
    \ann_{m-1}\{a,-1\}$, and we see that condition $(2)$ is
    indeed hereditary.
\end{proof}

\section{Examples of $H^k(E)$ Trivial for all $n<k$ and $H^n(E)$
Nontrivial, with Given Cohomological Dimension}\label{se:extriv}

We have shown in Theorem~\ref{th:trivher} that the property
$H^n(E)$ is a trivial $\Fp[G]$-module is hereditary. This result
leads naturally to the definition of an interesting invariant
$\ct(E/F) \in \{0\} \cup \N \cup\{\infty\}$:
\begin{multline*}
    \ct(E/F) \ = \sup \left\{ n\in\N \cup \{0\} \ \vert \
    H^n(E) \text{ is not a trivial } \Fp[G]\text{-module}\right\}.
\end{multline*}
As with $\cf$, we have chosen $\ct$ to indicate that after degree
$\ct(E/F)$, Galois cohomology consists of trivial
$\Fp[G]$-modules. Of course, if $H^n(E)$ is never trivial for
$n\ge 1$ then $\ct(E/F)=\infty$, and otherwise $\ct(E/F) \in
\N\cup \{0\}$. (Observe that since $H^0(E)\cong \Fp$ and there are
no nontrivial $G$-actions on $\Fp$, we always have that $H^0(E)$
is a trivial $\Fp[G]$-module. However, Theorem~\ref{th:trivher}
establishes the hereditary property only when $n > 0$.)

If $\ct(E/F)=n\in\N$, then by definition $H^m(E)$ is a trivial
$\Fp[G]$-module for all $m>n$. On the other hand, by the hereditary
property we also have that $H^k(E)$ is not trivial for all $1\le
k\leq \ct(E/F)$. A natural question arises: can we choose a suitable
field extension $E/F$ so that $\ct(E/F)$ is a given natural number
or $\infty$? We show that the answer is affirmative.  In fact, we
can arrange that both values $\ct(E/F)$ and $\cd(G_E)$ are any
natural numbers or $\infty$ and the absolute Galois group $G_F$ is a
pro-$p$-group modulo an obvious restriction, the inequality
described below.

Suppose that the absolute Galois group $G_E$ is a pro-$p$-group.
Then, observing that $\{0\}^G=\{0\}$, we have:
\begin{equation*}
    \ct(E/F) \leq \cd(G_E).
\end{equation*}

Our result is then the following.

\textit{Given $1\leq n\leq m\in\N\cup\{\infty\}$ and a prime $p$,
there exists a cyclic extension $E/F$ of degree $p$ with $\xi_p
\in F$ such that
\begin{enumerate}
    \item $G_E$ is a pro-$p$-group;
    \item $\ct(E/F)=n$; and
    \item $\cd(G_E)=m$.
\end{enumerate}
}

It is quite an interesting feature of our construction that it
parallels the construction made in the rather opposite free case
dealt with before. The only difference is the choice of $a$ in our
field extension of the form $E=F(\root{p}\of{a})$.

\subsection{The case $m\in\N$.}\

(1). Let $F:=F_{n,m}$ be a field of characteristic $0$ with $G_F
\cong \Z_p^n \star_{\text{pro-}p} \Z_p^m$ and $\xi_{p^2} \in F$,
given by Lemma~\ref{le:fmm}. Let
\begin{equation*}
    E=F(\root{p}\of{a})
\end{equation*}
where $a \in F^\times$ such that under the restriction map on $H^1$,
\begin{equation*}
    \res_{\star} (a) = (a)_1\oplus (a)_2, \qquad
    (a)_1 = 0, \ (a)_2 \neq 0.
\end{equation*}
Observe that there exists an $a$ with the required conditions
because by Lemma~\ref{le:fmlemma}, $H^1(F_m)\neq \{0\}$.

(2a). \textit{$H^n(E)$ is not trivial}. We claim that
\begin{equation*}
    \ann_{n} (a) \not\subset (a) \cup H^{n-1}(F).
\end{equation*}
By Lemma~\ref{le:fmlemma}, $H^n(F_n)$ contains a nontrivial
element $c$. Let $b\in H^n(F)$ such that
\begin{equation*}
    b_1=c\in H^n(F_n)\mbox{ and }b_2=0\in H^n(F_m).
\end{equation*}
Then $b\neq 0$ and since the cup-product commutes with $\res_
{\star}$,
\begin{equation*}
    \res_{\star} (a) \cup b = \big(0 \cup b_1\big) \oplus
    \big((a_1) \cup 0\big) = 0 \in H^{n+1} (F).
\end{equation*}
Therefore $b \in \ann_n (a)$.

Not let $f\in H^{n-1}(F)$ be arbitrary. Then
\begin{equation*}
    \left((a) \cup f \right)_1 = 0 \cup f_1 = 0
\end{equation*}
and therefore $b \notin (a) \cup H^{n-1}(F)$. Thus $\ann_n
(a) \not\subset (a) \cup H^{n-1}(F)$.

For the case $p>2$, Theorem~\ref{th:trivial}, part (2) implies
that $H^n(E)$ is not trivial.

In the case $p=2$ we have $(-1)=0$ since $\sqrt{-1}\in F^\times$.
Therefore $0=(a,-1)\in H^2(F)$. Thus $\ann_{n-1} (a,-1) = H^{n-1}
(F)$ and $(a) \cup \ann_{n-1} (a,-1) = (a) \cup H^{n-1}(F)$.
Hence by our claim above $\ann_n (a) \not\subset (a) \cup
\ann_{n-1} (a,-1)$, and we can again apply Theorem~\ref{th:trivial}
to conclude that $H^n(E)$ is not trivial.

(2b). \textit{$H^k(E)^G=H^k(E)$ for all $k\geq n+1$}.

Let $a_1\in F_m^\times$ satisfy $(a_1)=(a)_2$ and extend $\{(a_1)
\}$ to a basis $\{(a_1),\dots,(a_m)\}$ of $H^1(F_m)$. Recall that by
Lemma~\ref{le:fmlemma}, $H^k(F_m)$ is just the $k$th homogenous
summand of the exterior algebra over $\F_p$ generated by $H^1
(F_m)$. Using this fact and writing each element in $H^k(F_m)$ as a
sum of elements of the form
\begin{equation*}
    \left(a_{i_1}\right) \cup \cdots \cup \left(
    a_{i_k}\right), \ 1 \leq i_1 < i_2 < \cdots <
    i_k\leq m,
\end{equation*}
and also the fact that $H^k(F_n)=\{0\}$ we see that
\begin{equation*}
    \ann_k (a) = (a) \cup H^{k-1}(F).
\end{equation*}

Now again using Theorem~\ref{th:trivial} as in the case (2a), we
conclude that $H^k(E)^G=H^k(E)$.

(3). $\cd(G_E)=m$. Indeed $\cd(G_E)=\cd(G_F)$ by Serre's theorem.
(See \cite{S} and the discussion in section~4.1 which guarantees
that the hypothesis of Serre's theorem is valid.)

But from Lemma~\ref{le:fmlemma} and Lemma~\ref{le:fmm} we see that
    \begin{equation*}
        \cd(G_F)=\max\{\cd(G_{F_n}),\cd(G_{F_m})\}=m.
    \end{equation*}

Thus we see that in the case when $m<\infty$ we constructed a
cyclic field extension $E/F$ of degree $p$ with required properties
$(1),(2)$ and $(3)$.

\subsection{The case $m=\infty$.}

We first consider the subcase of this case when $n <\infty$. As in
section~\ref{ss:nminfty} set
    \begin{equation*}
        F_\infty:=\C((\Z_{(p)}^m)),\mbox{ where }
        m=\aleph_0.
    \end{equation*}
By Lemma~\ref{le:fmm} we see that there exists a field $F:=F_
{n,\infty}$ such that $G_F\cong G_{F_n} \star_{\text{pro-}p}
G_{F_\infty}$ and $\xi_{p^2} \in F$. Let $a \in F^\times$ such
that under the restriction map
    \begin{equation*}
        \res_{\star}: H^1(F)\to H^1(F_n)\oplus H^1(F_
        \infty)
    \end{equation*}
we have
    \begin{equation*}
        \res_{\star} (a) = 0 \oplus (a)_2, \ (a)_2 \ne 0.
    \end{equation*}
Then $\cd(F)=\infty$ and with the same argument as above we see
that $\ct(E/F)=n$.

Finally we consider the case $n=\infty=m$. Set again
$F_\infty:=\C((\Z_{(p)}^m))$, where $m=\aleph_0$ and $F=
F_{\infty,\infty}$. Also let $a\in F^\times$ such that
    \begin{equation*}
        \res_{\star} (a) = 0 \oplus (a)_2, (a)_2
        \ne 0.
    \end{equation*}
Then using the same argument as in (2b) we see that $\ct(F)=
\infty$.

Our construction is now completed.

\end{document}